\documentclass[10pt,twoside]{article}
\usepackage{amssymb}
\usepackage{amsmath}
\usepackage{Latex-document}

\newtheorem{sats}{Theorem}
\newtheorem{prop}{Proposition}
\newtheorem{lem}{Lemma}
\newtheorem{kor}{Corollary}
\newcommand{\banm}{\begin{anm}}
\newcommand{\eanm}{\end{anm}}

\newcommand{\tc}{\mathaccent"7017}
\newcommand{\Ho}{{\tc H}}

\title{\bf The Wiener Test for Higher Order\vskip -2mm Elliptic Equations\vskip 6mm}
\author{Vladimir Maz'ya\thanks{Link\"opings Universitet, Mathematiska
Institutionen, 58183 Link\"oping, Sweden. E-mail:\newline
vlmaz@mail.liu.se} \vspace*{-0.5cm}}
\date{\vspace{-8mm}}

\begin{document}

\maketitle

\thispagestyle{first} \setcounter{page}{189}

{\bf 1. Introduction.} Wiener's criterion for the regularity of a
boundary point with respect to the Dirichlet problem for the
Laplace equation [W] has been extended to various classes of
elliptic and parabolic partial differential equations. They
include linear divergence and nondivergence equations with
discontinuous coefficients, equations with degenerate quadratic
form, quasilinear and fully nonlinear equations, as well as
equations on Riemannian manifolds, graphs, groups, and metric
spaces (see [LSW], [FJK], [DMM], [LM], [KM], [MZ], [AH], [Lab],
[TW] to mention only a few). A common feature of these equations
is that all of them are of second order, and Wiener type
characterizations for higher order equations have been unknown so
far. Indeed, the increase of the order results in the loss of the
maximum principle, Harnack's inequality, barrier techniques, and
level truncation arguments, which are ingredients in different
proofs related to the Wiener test for the second order equations.

In the present work we extend Wiener's result to elliptic
differential operators $L(\partial)$ of order $2m$ in the
Euclidean space ${\bf R}^n$ with constant real coefficients
$$L(\partial)=(-1)^m\sum_{|\alpha|=|\beta|=m}a_{\alpha\beta}
 \partial^{\alpha+ \beta}.$$
We assume without loss of generality that $a_{\alpha\beta}=
a_{\beta\alpha}$ and $(-1)^mL(\xi)>0$ for all nonzero $\xi\in {\bf
R}^n$. In fact, the results can be extended to equations with
variable (for example, H\"older continuous) coefficients in
divergence form but we leave aside this generalization to make
exposition more lucid.

We use the notation $\partial$ for the gradient
$(\partial_{x_1},\ldots ,\partial_{x_n})$, where $\partial_{x_k}$
is the partial derivative with respect to $x_k$. By $\Omega$ we
denote an open set in ${\bf R}^n$ and by $B_{\rho}(y)$ the ball
$\{x\in {\bf R}^n: |x-y|<\rho\}$, where $y\in {\bf R}^n$. We write
$B_{\rho}$ instead of $B_{\rho}(O)$.

Consider the Dirichlet problem
\begin{equation} \label{dir}
  L(\partial)u=f, \;\; f\in C^{\infty}_0(\Omega), \;\; u\in
  \Ho^m(\Omega),
\end{equation}
where we use the standard notation $C^{\infty}_0(\Omega)$ for the
space of infinitely differentiable functions in ${\bf R}^n$ with
compact support in $\Omega$ as well as $\Ho^m(\Omega)$ for the
completion of $C^{\infty}_0(\Omega)$ in the energy norm.

We call the point $O\in \partial\Omega$ regular with respect to
$L(\partial)$ if for any $f\in C^{\infty}_0(\Omega)$ the solution
of (\ref{dir}) satisfies
\begin{equation} \label{c69}
  \mathop{\hbox {lim}}_{\Omega\ni x\to O}u(x)=0.
\end{equation}

For $n=2,3, \ldots , 2m-1$ the regularity is a consequence of the
Sobolev imbedding theorem. Therefore, we suppose that $n\geq 2m$.
In the case $m=1$ the above definition of regularity is equivalent
to that given by Wiener.

The following result coincides with Wiener's criterion in the case
$n=2$ and $m=1$.

\begin{sats}\label{Th1} Let $2m=n$. Then $O$ is regular with respect
to $L(\partial)$ if and only if
\begin{equation} \label{c71}
  \int _0^1C_{2m}(B_{\rho}\backslash \Omega)\rho^{-1}d\rho=\infty.
\end{equation}
\end{sats}

Here and elsewhere $C_{2m}$ is the potential-theoretic Bessel
capacity of order $2m$ (see [AHed]). If $n=2m$ and $O$ belongs to
a continuum contained in the complement of $\Omega$, condition
(\ref{c71}) holds.

The case $n>2m$ is more delicate because no result of Wiener's
type is valid for all operators $L(\partial)$ (see [MN]). To be
more precise, even the vertex of a cone can be irregular with
respect to $L(\partial)$ if the fundamental solution of
$L(\partial)$:
\begin{equation} \label{c70}
  F(x)=F(x/|x|)|x|^{2m-n}, \;\; x\in {\bf R}^n\backslash O,
\end{equation}
changes sign. Examples of operators $L(\partial)$ with this
property were given in [MN] and [D]. For instance, according to
[MN] the vertex of a sufficiently thin 8-dimensional cone $K$ is
irregular with respect to the operator
\[ L(\partial)u:= 10\partial_{x_8}^4u + \Delta^2u,
~u\in \Ho^2({\bf R}^8\backslash K). \]

In the sequel, Wiener's type characterization of regularity for
$n>2m$ is given for a subclass of the operators $L(\partial)$
called {\it positive with the weight} $F$. This means that for all
real-valued $u\in C_0^{\infty}({\bf R}^n\backslash O)$,
\begin{equation} \label{c4}
  \int_{{\bf R}^n}L(\partial)u(x)\cdot u(x)F(x)\,dx\geq c\;
  \sum_{k=1}^m\int_{{\bf R}^n}|\nabla _ku(x)|^2|x|^{2k-n}dx,
\end{equation}
where $\nabla _k$ is the gradient of order $k$, i.e.
$\nabla_k=\{\partial^{\alpha}\}$ with $|\alpha|=k$.

The positivity of the left-hand side in (\ref{c4}) is equivalent
to the inequality
$$\int_{{\bf R}^n}\int_{{\bf R}^n} \frac{L(\xi)+L(\eta)}{L(\xi-\eta)}
f(\xi)f(\eta)d\xi d\eta>0$$ for all non-zero $f\in
C^{\infty}_0({\bf R}^n)$.

\begin{sats}\label{Th2} Let $n>2m$ and let $L(\partial)$ be
positive with weight $F$. Then $O$ is regular with respect to
$L(\partial)$ if and only if
\begin{equation} \label{c72}
  \int _0^1C_{2m}(B_{\rho}\backslash \Omega)\rho^{2m-n-1}d\rho=
  \infty.
\end{equation}
\end{sats}

Note that in direct analogy with the case of the Laplacian we
could say, in Theorems \ref{Th1} and \ref{Th2}, that $O$ is
irregular with respect to $L(\partial)$ if and only if the set
${\bf R}^n\backslash\Omega$ is $2m$-thin in the sense of linear
potential theory [L], [AHed].

Let, for example, the exterior of $\Omega$ contain the region
$$\{\, x: 0<x_n<1, \quad \left(x_1^2+ \ldots +
x_{n-1}^2\right)^{1/2}<f(x_n)\,\},$$
where $f$ is an increasing
function such that $f(0)=f'(0)=0$. Then the point $O$ satisfies
(\ref{c72}) if and only if
$$ \int_0^1|\log f(\tau)|^{-1}\tau^{-1}d\tau=\infty \quad {\rm for}
\;\; n=2m+1$$
and
$$ \int_0^1f(\tau)\tau^{2m-n}d\tau=\infty \quad \quad \quad
{\rm for} \;\; n\geq 2m+2.$$

Since, obviously, the operator $L(\partial)$ of the second order
is positive with the weight $F$, Wiener's result for $n>2$ is
contained in Theorem 2.

We note that the pointwise positivity of $F$ follows from
(\ref{c4}), but the converse is not true. In particular, the
$m$-harmonic operator with $2m<n$ satisfies {\rm (\ref{c4})} if
and only if $n=5,6,7$ for $m=2$ and $n=2m+1$, $2m+2$ for $m>2$
(see [M2], where the proof of sufficiency of {\rm (\ref{c72})} is
given for $(-\Delta)^m$ with $m$ and $n$ as above, and also [E]
dealing with the sufficiency for noninteger powers of the
Laplacian in the intervals $(0,1)$ and $[n/2-1,n/2)$).

We state some auxiliary assertions of independent interest which
concern the so called $L$-capacitary potential $U_K$ of the
compact set $K\subset {\bf R}^n$, $n>2m$, i.e. the solution of the
variational problem
$$ {\rm inf}\{\int_{{\bf R}^n}L(\partial)u\cdot u\,dx:
\;\; u\in C_0^{\infty}({\bf R}^n), u=1 \; {\rm in} \; {\rm
vicinity} \; {\rm of} \;K\}.$$

These assertions are used in the proof of necessity in Theorem 2.

By the $m$-harmonic capacity ${\rm cap}_m(K)$ of a compact set $K$
we mean
\begin{equation} \label{w51}
  {\rm inf}\Bigl\{\sum_{|\alpha|=m} \frac{m!}{\alpha!}
  ||\partial^{\alpha} u||^2_{L_2({\bf R}^n)}: \quad u\in
  C^{\infty}_0(\Omega), \; u=1 \;\; {\rm in} \; {\rm vicinity}
  \; {\rm of} \; K\Bigr\}.
\end{equation}

\begin{lem}\label{L1} Let $\Omega={\bf R}^n$, $2m<n$. For all $y
\in {\bf R}^n\backslash K$
\begin{eqnarray} \label{a3}
  &&U_K(y)=2^{-1}U_K(y)^2 \nonumber\\
  &&+\int_{{\bf R}^n}\; \sum_{m\geq
  j\geq1}\sum_{|\mu|=|\nu|=j}\partial^{\mu}U_K(x)\!\cdot\!
  \partial^{\nu}U_K(x)\!\cdot \!{\cal
  P}_{\mu\nu}(\partial)F(x-y)\,dx,
\end{eqnarray}
where ${\cal P}_{\mu\nu}(\zeta)$ are homogeneous polynomials of
degree $2(m-j)$, ${\cal P}_{\mu\nu}={\cal P}_{\nu\mu}$ {\it and}
${\cal P}_{\alpha\beta}(\zeta)=a_{\alpha\beta}$ for
$|\alpha|=|\beta|=m$.
\end{lem}

\begin{kor}\label{P1} Let $\Omega={\bf R}^n$ {\it and} $2m<n$.
For all $y\in {\bf R}^n\backslash K$ there holds the estimate
\begin{equation} \label{a5}
  |\nabla _jU_K(y)|\leq c_j\; {\rm dist}(y,K)^{2m-n-j}\; {\rm
  cap}_mK,
\end{equation}
where $j=0,1,2, \ldots$ and $c_j$ does not depend on $K$ and $y$.
\end{kor}

By ${\cal M}$ we denote the Hardy-Littlewood maximal operator.

\begin{kor}\label{P2} Let $2m<n$ and let $0<\theta<1$. Also
let $K$ be a compact subset of $\overline{B_{\rho}}\backslash
B_{\theta\rho}$. Then the $L$-capacitary potential $U_K$ satisfies
\begin{equation} \label{a8}
  {\cal M}\nabla _lU_K(0)\leq c_{\theta}\; \rho^{2m-l-n}\;
  {\rm cap}_mK,
\end{equation}
where $l=0, 1, \ldots , m$ and $c_{\theta}$ does not depend on $K$
and $\rho$.
\end{kor}

Let $L(\partial)$ be positive with the weight $F$. Then identity
{\rm (\ref{a3})} implies that the $L$-capacitary potential of a
compact set $K$ with positive $m$-harmonic capacity satisfies
\begin{equation} \label{c7}
  0<U_K(x)<2 \;\;\; {\rm on} \;\; {\bf R}^n\backslash K.
\end{equation}
In general, the bound $2$ in {\rm (\ref{c7})} cannot be replaced
by $1$.

\begin{prop}\label{P6} If $L=\Delta ^{2m}$, then there exists a
compact set $K$ such that
$$(U_K-1)\big |_{{\bf R}^n\backslash K}$$
changes sign in any neighbourhood of a point of $K$.
\end{prop}

We give a lower pointwise estimate for $U_K$ stated in terms of
capacity (compare with the upper estimate {\rm (\ref{a5}))}.

\begin{prop}\label{P7} Let $n>2m$ and let $L(\partial)$ be
positive with the weight $F$. If $K$ is a compact subset of $B_d$
and $y\in {\bf R}^n\backslash K$, then
$$U_K(y)\geq c\; (|y|+d)^{2m-n}\; {\rm cap} _mK.$$
\end{prop}

Sufficiency in Theorem 2 follows from the next assertion which is
of interest in itself.

\begin{lem}\label{L5} Let $2m < n$ and let $L(\partial)$ be
positive with the weight $F$. Also let $u \in \Ho^{m}(\Omega)$
satisfy $L(\partial) u = 0$ on $\Omega \cap B_{2R}$. Then, for all
$\rho \in (0,R)$,
\begin{eqnarray} \label{n4}
  &&\sup\{|u(p)|^2 : p \in \Omega \cap B_\rho \} +
  \int_{\Omega \cap B_\rho} \sum^m_{k=1} \frac{|\nabla_ku(x)|^2}{|x|^{n-2k}}\,dx
  \nonumber\\
  &&\leq c_1M_R(u)\exp \Bigl(-c_2\int^R_\rho {\rm
  cap}_m(\bar{B}_\tau\setminus \Omega) \frac{d\tau}{\tau^{n-2m+1}} \Bigr),
\end{eqnarray}
where $c_1$ and $c_2$ are positive constants, and
$$ M_R(u)=R^{-n}\int_{\Omega\cap(B_{2R}\backslash B_R)}|u(x)|^2dx. $$
\end{lem}

The present work gives answers to some questions posed in [M2]. I
present several simply formulated unsolved problems.

{\it 1. Is it possible to replace the positivity of} $L(\partial)$
{\it with the weight} $F(x)$ {\it by the positivity of} $F(x)$
{\it in Theorem 2}?

A particular case of this problem is the following one.

{\it 2. Does Theorem 2 hold for the operator} $(-\Delta)^m$, {\it
where}
\begin{equation}\label{n1}
  n\geq 8, \;\; m=2 \quad {\rm and} \quad n\geq 2m+3, \; m>2 \;?
\end{equation}

The next problem concerns Green's function $G_m$ of the Dirichlet
problem for $(-\Delta)^m$ in an arbitrary domain $\Omega$.

{\it 3. Prove or disprove the estimate}
\begin{equation}\label{n2}
  |G_m(x,y)|\leq \frac{c(m,n)}{|x-y|^{n-2m}},
\end{equation}
{\it where} $c(m,n)$ {\it is independent of} $\Omega$ {\it and}
$m$ {\it and} $n$ {\it are the same as in} (\ref{n1}).

For $n=5, 6, 7$, $m=2$ and $n=2m+1$, $2m+2$, $m>2$ estimate
(\ref{n1}) was proved in [M3]. In the sequel, by $u$ we denote a
solution in $\Ho^{m}(\Omega)$ of the equation
\begin{equation} \label{15}
  (-\Delta)^m u = f \quad {\rm in} \quad \Omega.
\end{equation}
Clearly, (\ref{n2}) leads to the following estimate of the maximum
modulus of $u$
$$\|u\|_{L_{\infty}}\leq c(m,n,{\rm mes}_n\Omega)\|f\|_{L_p(\Omega)},$$
where $p>n/2m$. However, the validity of this estimate for the
same $n$ and $m$ as in (\ref{n1}) is an open problem. Moreover,
the following questions arise.

{\it 4. Let} $m=2$, $n\geq 8$, {\it and let} $\Omega$ {\it be an
arbitrary bounded domain.} {\it Is} $u$ {\it uniformly bounded in}
$\Omega$ {\it for any} $f\in C_0^{\infty}(\Omega)$?

{\it 5. Let} $m>2$ {\it and} $n\geq 2m+3$. {\it Also, let}
$\partial\Omega$ {\it have a conic singularity.} {\it Is} $u$ {\it
uniformly bounded in} $\Omega$ {\it for any} $f\in
C_0^{\infty}(\Omega)$?

For $m=2$, the affirmative answer to the last question is given
in [MP].

I formulate two related open problems.

{\it 6. Let} $m=2$ {\it and} $n=2$. {\it Is } $u$ {\it Lipschitz
up to the boundary of an arbitrary bounded domain, for any} $f\in
C_0^\infty (\Omega)$ ?

{\it 7. Let } $m=2$ and $n\geq 3$. {\it Does} $u$ {\it belong to
the class} $C^{1,1}(\Omega)$ {\it for any} $f\in
C_0^\infty(\Omega)$ {\it if} $\Omega$ {\it is convex}?

According to [KoM], the last is true in the two-dimensional case.

I conclude with the following variant of the Phragm\'en-Lindel\"of
principle (see [M3]).

\begin{prop}\label{P3} Let either $n=5,6,7$, $m=2$ or $n=2m+1$,
$2m+2$, $m>2$. Further, let $\eta u\in \Ho^{m}(\Omega)$ for all
$\eta\in C^{\infty}({\bf R}^n)$, $\eta=0$ near $O$. If
$$\Delta ^m u=0 \quad {\rm on } \; \; \Omega\cap B_1,$$
then either $u\in \Ho^{m}(\Omega)$ and
$$ \mathop{\operatorname{lim\ sup}}_{\rho\to 0}\sup_{B_{\rho}\cap\Omega}
|u(x)|\exp \Bigl(c\int^1_\rho {\rm cap}_m(\bar{B}_\rho\setminus
\Omega) \frac{d\rho}{\rho} \Bigr)<\infty$$ or
$$\mathop{\operatorname{lim\ inf}}_{\rho\to 0}\rho^{n-2m}M_{\rho}(u)\exp
\Bigl(-c\int^1_\rho {\rm cap}_m(\bar{B}_\rho\setminus \Omega)
\frac{d\rho}{\rho} \Bigr)>0. $$
\end{prop}

It would be interesting to extend this assertion to other values
of $n$ and $m$.

\label{lastpage}

\end{document}